\documentclass[11pt]{article}

\usepackage[dcucite]{harvard}
\renewcommand{\harvardyearleft}{}
\renewcommand{\harvardyearright}{}

\usepackage{amsmath}
\usepackage{amsthm}
\usepackage{amscd}
\usepackage{amsbsy}
\usepackage{amssymb}
\usepackage{amsmath}
\usepackage{amsfonts}
\usepackage{amsopn}
\usepackage{a4wide}

\usepackage[latin1]{inputenc}

\theoremstyle{plain}
\newtheorem{theorem}{Theorem}[section]
\newtheorem{corollary}[theorem]{Corollary}

\theoremstyle{definition}

\theoremstyle{remark}
\newtheorem{remark}{Remark}[section]

\begin{document}

\title{L'H\^{o}pital-Type Rules for Monotonicity\\
with Application to Quantum Calculus\thanks{Submitted 25-Jul-2010; accepted 11-Nov-2010. 
Ref: Int. J. Math. Comput. 10 (2011), no. M11, 99--106.}}

\author{\textbf{Nat\'{a}lia Martins$^1$ and Delfim F. M. Torres$^2$}}

\date{$^1$Department of Mathematics\\
University of Aveiro\\
3810-193 Aveiro, Portugal\\
natalia@ua.pt\\[0.3cm]
$^2$Department of Mathematics\\
University of Aveiro\\
3810-193 Aveiro, Portugal\\
delfim@ua.pt}

\maketitle


\begin{abstract}
\noindent \emph{We prove new l'H\^{o}pital rules for
monotonicity valid on an arbitrary time scale. Both delta and nabla monotonic
l'H\^{o}pital rules are obtained. As an example of application,
we give some new upper and lower bounds for the exponential
function of quantum calculus restricted to the $q$-scale.}

\medskip

\noindent\textbf{Keywords:} l'H\^{o}pital-type rules,
monotonic functions, $q$-calculus.

\medskip

\noindent\textbf{2010 Mathematics Subject Classification:} 26A48, 39A13.

\end{abstract}


\section{Introduction}

L'H\^{o}pital-type rules for monotonicity have an important
role in mathematics \cite{Anderson:et:al,Boas,Pinelis:2006,Pinelis:2008},
with numerous useful applications found in mathematical inequalities,
statistics, probability, differential geometry, approximation and
information theories, and mathematical physics
\cite{Andras,Baricz,Baricz2010,Pan,Pinelis:2004,Pinelis:2007,Zhu}.
Recently, Wu and Debnath proved the following
l'H\^{o}pital-type rules for monotonicity:

\begin{theorem}[\cite{S.-Wu}]
\label{Hospital rule}
Let $f$ and $g$ be differentiable functions on $]a,b[$. Suppose
that either $g^\prime >0$ everywhere on $]a,b[$ or $g^\prime <0$
everywhere on $]a,b[$.

If $\displaystyle\frac{f^\prime}{g^\prime}$ is increasing (resp.
decreasing) on $]a,b[$ and $f(a^+)$ and $g(a^+)$ exist, then the
function
$$\frac{f(x)-f(a^+)}{g(x)-g(a^+)}$$
is also increasing (resp. decreasing) on $]a,b[$.

\clearpage

If $\displaystyle\frac{f^\prime}{g^\prime}$ is increasing (resp.
decreasing) on $]a,b[$ and $f(b^-)$ and $g(b^-)$ exist, then the
function
$$\frac{f(x)-f(b^-)}{g(x)-g(b^-)}$$
is also increasing (resp. decreasing) on $]a,b[$.
\end{theorem}

In this note we generalize Theorem~\ref{Hospital rule}
by proving l'H\^{o}pital-type rules for monotonicity
on an arbitrary time scale (Section~\ref{sec:MR}).
Our results seem to be new even in the discrete-time case.
As an application, we prove new inequalities on quantum calculus:
we give new upper and lower bounds for the exponential
function of quantum calculus when restricted to the $q$-scale
(Section~\ref{sec:Appl}).

The theory of time scales was introduced in 1988 by Aulbach
and Hilger in order to unify continuous and
discrete-time theories \cite{AH88}. It has found applications in several
different fields that require simultaneous modeling of discrete
and continuous data \cite{Guseinov-2002,S.-Wu}, and is now
a subject of strong current research (see
\cite{Ricardo+Delfim,Zbig+Delfim,Agnieszka+Delfim,Natalia+Delfim,Dorota+Delfim,Rachid+Rui+Delfim}
and references therein). We claim that time scale theory is useful with respect
to L'H\^{o}pital-type rules for monotonicity, avoiding repetition
of results in the continuous and discrete cases \cite{Pinelis:2007,Pinelis:2008}.
We trust that the present study will mark the beginning
of the study of L'H\^{o}pital-type rules on time scales
and its applications to dynamic inequalities \cite{GronTS,Rachid+Rui+Delfim,SA:T},
the calculus of variations \cite{Bohner:Rui:Delfim,Rui:Rachid:Delfim},
and quantum calculus \cite{Cresson:Gastao:Delfim,Victor-Kac-2002}.


\section{L'H\^{o}pital rules for monotonicity}
\label{sec:MR}

For definitions, notations, and results concerning
the theory of delta and nabla calculus
on time scales we refer the reader to the books
\cite{Bohner-Peterson-2001,Bohner-Peterson-2003}.
Similarly to \cite{S.-Wu}, we use $f(a^+)$ to denote
$\lim_{x\rightarrow a^+} f(x)$ and $f(b^-)$ to denote
$\lim_{x\rightarrow b^-} f(x)$. All the intervals in this work
are time scale intervals.

\begin{theorem}[Delta-monotonic l'H\^{o}pital rules]
\label{hospital thm delta}
Let $\mathbb{T}$ be a time scale with $a,b \in \mathbb{T}$, $a<b$,
and $f$ and $g$ two $\Delta$-differentiable functions on $]a,b[$.
Suppose that either $g^\Delta >0$ everywhere on $]a,\rho(b)[$ or
$g^\Delta <0$ everywhere on $]a,\rho(b)[$.
\begin{enumerate}
\item  If $\displaystyle\frac{f^\Delta}{g^\Delta}$ is increasing
(resp. decreasing) on $]a,\rho(b)[$ and $f(a^+)$ and $g(a^+)$
exist (finite), then the function
$$\frac{f(x)-f(a^+)}{g(x)-g(a^+)}$$
is also increasing (resp. decreasing) on $]a,\rho(b)[$.
\item If $\displaystyle\frac{f^\Delta}{g^\Delta}$ is increasing
(resp. decreasing) on $]a,\rho(b)[$ and $f(b^-)$ and $g(b^-)$
exist (finite), then the function
$$\frac{f(x)-f(b^-)}{g(x)-g(b^-)}$$
is also increasing (resp. decreasing) on $]a,\rho(b)[$.
\end{enumerate}
\end{theorem}

\begin{proof}
We will prove the first assertion.
Since $f(a^+)$ and $g(a^+)$ are finite numbers,
we can define the following functions on $[a,b[$:
$$
F(x)=  \left\{ \begin{array}{lcl}
f(x)& \mbox{if} & x \in ]a,b[\\
f(a^+)& \mbox{if} & x=a \end{array} \right.
$$
and
$$
G(x)=  \left\{ \begin{array}{lcl}
g(x)& \mbox{if} & x \in ]a,b[\\
g(a^+)& \mbox{if} & x=a \, .
\end{array} \right.
$$
Clearly, for any $x \in ]a,b[$, $F$ and $G$ are continuous on
$[a,x]$ and $\Delta$-differentiable  on $[a,x[$. Note also that
$G(x)\neq G(a)$ for $x \ne a$.

By the Cauchy mean value theorem with delta derivatives
\cite{Guseinov-2002} we conclude that there exist
$c_1,c_2 \in [a,x[$ such that
$$
\frac{F^\Delta(c_1)}{G^\Delta(c_1)} \leq
\frac{F(x)-F(a)}{G(x)-G(a)} \leq
\frac{F^\Delta(c_2)}{G^\Delta(c_2)}.
$$
For each  $x \in ]a,b[$ define
$$
H(x)=\frac{F(x)-F(a)}{G(x)-G(a)}.
$$
Note that, for each $x \in ]a,\rho(b)[$,
$$
\begin{array}{rcl}
\displaystyle H^\Delta(x)& = & \displaystyle \frac{F^\Delta(x)
(G(x)-G(a))- G^\Delta(x)
(F(x)-F(a))}{(G(x)-G(a))(G^{\sigma}(x)-G(a))}\\
& &\\
& = & \displaystyle
 \frac{G^\Delta(x)}{G^{\sigma}(x)-G(a)}\left(\frac{F^\Delta(x)}{G^\Delta(x)
}-\frac{F(x)-F(a)}{G(x)-G(a)}\right).
\end{array}
$$

Case I. Suppose that $\displaystyle\frac{f^\Delta}{g^\Delta}$ is
increasing on $]a,\rho(b)[$. For each $x \in ]a,\rho(b)[$ there
exists $c_2 \in [a,x[$ such that
$$
\begin{array}{rcl}
\displaystyle H^\Delta(x) &\geq& \displaystyle
\frac{G^\Delta(x)}{G^{\sigma}(x)-G(a)}\left(\frac{F^\Delta(x)}{G^\Delta(x)
}-\frac{F^\Delta(c_2)}{G^\Delta(c_2)}\right)\\
&=& \displaystyle
\frac{g^\Delta(x)}{g^{\sigma}(x)-g(a^+)}\left(\frac{f^\Delta(x)}{g^\Delta(x)
}-\frac{f^\Delta(c_2)}{g^\Delta(c_2)}\right)>0.
\end{array}
$$
Hence, $H$ is increasing on $]a,\rho(b)[$ \cite{Guseinov-2002}. Since
$$H(x)=\frac{f(x)-f(a^+)}{g(x)-g(a^+)}$$ we can conclude that
$\displaystyle\frac{f(x)-f(a^+)}{g(x)-g(a^+)}$ is increasing on
$]a,\rho(b)[$.

Case II. Suppose that $\displaystyle\frac{f^\Delta}{g^\Delta}$ is
decreasing on $]a,\rho(b)[$. For each $x \in ]a,\rho(b)[$ there
exists $c_1 \in [a,x[$ such that
$$
\begin{array}{rcl}
\displaystyle H^\Delta(x) &\leq& \displaystyle
\frac{G^\Delta(x)}{G^{\sigma}(x)-G(a)}\left(\frac{F^\Delta(x)}{G^\Delta(x)
}-\frac{F^\Delta(c_1)}{G^\Delta(c_1)}\right)\\
&=& \displaystyle
\frac{g^\Delta(x)}{g^{\sigma}(x)-g(a^+)}\left(\frac{f^\Delta(x)}{g^\Delta(x)
}-\frac{f^\Delta(c_1)}{g^\Delta(c_1)}\right)<0.
\end{array}
$$
Hence, $H$ is decreasing on $]a,\rho(b)[$.  This proves that
$\displaystyle\frac{f(x)-f(a^+)}{g(x)-g(a^+)}$ is decreasing on
$]a,\rho(b)[$.

The second assertion is proved in a similar way.
\end{proof}

From Theorem~\ref{hospital thm delta} we can deduce the following
corollaries.

\begin{corollary}
\label{Corollary 1}
Let $\mathbb{T}$ be a time scale with $a,b \in \mathbb{T}$, $a<b$,
and $f,g:[a,b[ \rightarrow\mathbb{R}$
be two continuous functions which are $\Delta$-differentiable on $]a,b[$.
Suppose that either $g^\Delta >0$ everywhere on $]a,\rho(b)[$ or
$g^\Delta <0$ everywhere on $]a,\rho(b)[$.

If $\displaystyle\frac{f^\Delta}{g^\Delta}$ is increasing (resp.
decreasing) on $]a,\rho(b)[$, then the function
$$\frac{f(x)-f(a)}{g(x)-g(a)}$$
is also increasing (resp. decreasing) on $]a,\rho(b)[$.
\end{corollary}

\begin{corollary}
\label{Corollary 2}
Let $\mathbb{T}$ be a time scale with $a,b \in \mathbb{T}$, $a<b$,
and $f,g : \ ]a,b] \rightarrow\mathbb{R}$
be two continuous functions which are $\Delta$-differentiable on $]a,b[$.
Suppose that either $g^\Delta >0$ everywhere on $]a,\rho(b)[$ or
$g^\Delta <0$ everywhere on $]a,\rho(b)[$.

If $\displaystyle\frac{f^\Delta}{g^\Delta}$ is increasing (resp.
decreasing) on $]a,\rho(b)[$, then the function
$$\frac{f(x)-f(b)}{g(x)-g(b)}$$
is also increasing (resp. decreasing) on $]a,\rho(b)[$.
\end{corollary}

\begin{corollary}
\label{Corollary 3}
Let $\mathbb{T}$ be a time scale with $a,b \in \mathbb{T}$, $a<b$,
and $f,g : \  ]a,b[ \rightarrow\mathbb{R}$ be two
$\Delta$-differentiable functions on $]a,b[$.
Suppose that either $g^\Delta >0$ everywhere on $]a,\rho(b)[$ or
$g^\Delta <0$ everywhere on $]a,\rho(b)[$. Suppose also that
$f(a^+)=g(a^+)=0$ or $f(b^-)=g(b^-)=0$.

If $\displaystyle\frac{f^\Delta}{g^\Delta}$ is increasing (resp.
decreasing) on $]a,\rho(b)[$, then $\displaystyle\frac{f}{g}$ is
also increasing (resp. decreasing) on $]a,\rho(b)[$.
\end{corollary}

\begin{remark}
Theorem \ref{hospital thm delta} and Corollaries
\ref{Corollary 1}, \ref{Corollary 2} and \ref{Corollary 3} hold
true if the terms ``increasing'' and ``decreasing'' are replaced
everywhere by ``non-decreasing'' and ``non-increasing'', respectively.
\end{remark}

Using the recent duality theory \cite{Caputo,BD,MyID:177},
one can easily obtain the corresponding nabla result
for Theorem~\ref{hospital thm delta}:

\begin{theorem}[Nabla-monotonic l'H\^{o}pital rules]
\label{hospital thm nabla}
Let $\mathbb{T}$ be a time scale with $a,b \in \mathbb{T}$, $a<b$,
and $f$ and $g$ be $\nabla$-differentiable functions on $]a,b[$.
Suppose that either $g^\nabla >0$ everywhere on $]\sigma(a),b[$ or
$g^\nabla <0$ everywhere on $]\sigma(a),b[$.
\begin{enumerate}
\item  If $\displaystyle\frac{f^\nabla}{g^\nabla}$ is increasing
(resp. decreasing) on $]\sigma(a),b[$ and $f(a^+)$ and $g(a^+)$
exist (finite), then the function
$$\frac{f(x)-f(a^+)}{g(x)-g(a^+)}$$
is also increasing (resp. decreasing) on $]\sigma(a),b[$.

\item If $\displaystyle\frac{f^\nabla}{g^\nabla}$ is increasing
(resp. decreasing) on $]\sigma(a),b[$ and $f(b^-)$ and $g(b^-)$
exist (finite), then the function
$$\frac{f(x)-f(b^-)}{g(x)-g(b^-)}$$
is also increasing (resp. decreasing) on $]\sigma(a),b[$.
\end{enumerate}
\end{theorem}

In the case $\mathbb{T}=\mathbb{R}$
Theorems~\ref{hospital thm delta} and \ref{hospital thm nabla}
coincide and reduce to Theorem~\ref{Hospital rule}.
The l'H\^{o}pital rules of Theorems~\ref{hospital thm delta}
and \ref{hospital thm nabla}
seem to be new even for $\mathbb{T}=\mathbb{Z}$.
In the next section we give an application of our results
for $\mathbb{T}= \overline{q^\mathbb{Z}}=q^\mathbb{Z}\cup \{0\}$.


\section{An application to quantum inequalities}
\label{sec:Appl}

In this section we assume familiarity with the definitions
and results from the $q$-calculus (\textrm{cf.}, \textrm{e.g.},
the book \cite{Victor-Kac-2002}).
Let $q$ be a real number, $0<q<1$.
The \emph{$q$-derivative} $D_q f$
of a function $f$ is defined by
\begin{equation*}
D_q f (x)=
\begin{cases}
\displaystyle \frac{f(qx)-f(x)}{(q-1)x} & \text{ if } \ x\neq 0\\
f^\prime(x) & \text{ if } \ x = 0,
\end{cases}
\end{equation*}
provided $f^\prime(0)$ exists.
We use the following standard notations of $q$-calculus:

\begin{itemize}
\item for any real number $\alpha$,
$\displaystyle [\alpha]:= \displaystyle \frac{q^{\alpha}-1}{q-1}$;

\item  $[n]!:=
\begin{cases}
1 & \text{ if } n=0,\\
[n][n-1] \cdots [2][1] & \text{ if } n \in \mathbb{N};
\end{cases}$

\item $(x-a)^n_q:=
\begin{cases}
1 & \text{ if } n=0,\\
(x-a)(x-qa)\cdots(x-q^{n-1}a)  & \text{ if }  n \in \mathbb{N}.
\end{cases}$
\end{itemize}

The function $e^x_q$ defined by
\begin{equation}
\label{eq:q:exp}
e^x_q=\sum_{k=0}^{+\infty}\frac{x^k}{[k]!}
\end{equation}
is called the \emph{$q$-exponential function}.
The basic properties of the $q$-exponential function
can be found, \textrm{e.g.},
in \cite{Victor-Kac-2002,A.-Lavagno-2007}.
Here we obtain new upper and lower bounds for \eqref{eq:q:exp}.

\begin{theorem}
Let $0<q<1$, $a,b \in q^\mathbb{Z}$, $a<b$, and $n \in
\mathbb{N}$. Then, for any $x \in [q^{-1} a, b]$,
the following inequalities hold:
\begin{gather*}
e^x_q \geq \displaystyle \sum_{k=0}^{n-1} \frac{e_q^{a}}{[k]!}(x-a)^{k}_{q} +
 \frac{1}{(q^{-1}a-a)^n_q} \left( e^{q^{-1}a}_q - \displaystyle \sum_{k=0}^{n-1}
\frac{e_q^{a}}{[k]!}(q^{-1}a-a)^{k}_{q} \right) (x-a)_q^n,\\
e^x_q \leq \displaystyle \sum_{k=0}^{n-1}
\frac{e_q^{a}}{[k]!}(x-a)^{k}_{q} + \frac{1}{(b-a)^n_q} \left(
e^{b}_q - \displaystyle \sum_{k=0}^{n-1}
\frac{e_q^{a}}{[k]!}(b-a)^{k}_{q} \right) (x-a)_q^n.
\end{gather*}
\end{theorem}

\begin{proof}
We begin by noting that if we consider
$$
\mathbb{T}=\overline{q^\mathbb{Z}}=q^\mathbb{Z}\cup \{0\},
$$
then, for any function $f:\mathbb{T}\rightarrow\mathbb{R}$,
one obtains
$$
f^{\nabla}(t) = D_q f(t), \quad \forall t \in \mathbb{T},
$$
and $\sigma(t) = q^{-1} t$. The result follows as a corollary of
Theorem~\ref{hospital thm nabla}. For $x \in [a,b]$ define
$$
f(x)=e^x_q - \sum_{k=0}^{n-1}
\frac{e_q^{a}}{[k]!}(x-a)^{k}_{q}
$$
and
$$
g(x)=(x-a)^n_q.
$$
Note that
$$
D_q^{k} f(a)= f^{\nabla^{k}}(a)=0 \ \ \ \mbox{and} \ \ \
D_q^{k} g(a)=g^{\nabla^{k}}(a)=0 \ \ \ \mbox{for all} \ \
k=0,1,2,\ldots,n-1,
$$
where $D^k_q f$ is defined by $D^0_q f=f$
and $D_q^k f= D_q(D_q^{k-1})f$, $k=1,2,3,\ldots$
Note also that
$$
D_q^{n} f(x)=f^{\nabla^{n}}(x)= e^x_q \ \ \
\mbox{and} \ \ \  D_q^{n} g(x)=g^{\nabla^{n}}(x)= [n]!
$$
Thus, we have:
\begin{equation}
\label{formulas}
\begin{array}{rcl}
\displaystyle \frac{f(x)}{g(x)}& = & \displaystyle \frac{f(x)-f(a)}{g(x)-g(a)}\\
&&\\ \displaystyle \frac{f^{\nabla}(x)}{g^{\nabla}(x)} & = &
\displaystyle \frac{f^{\nabla}(x)-f^{\nabla}(a)}{g^{\nabla}(x)-g^{\nabla}(a)}\\
\displaystyle
&\vdots& \\
\displaystyle \frac{f^{\nabla^{n-1}}(x)}{g^{\nabla^{n-1}}(x)}& = &
\displaystyle \frac{f^{\nabla^{n-1}}(x)
-f^{\nabla^{n-1}}(a)}{g^{\nabla^{n-1}}(x)-g^{\nabla^{n-1}}(a)}\\
&&\\
\displaystyle \frac{f^{\nabla^{n}}(x)}{g^{\nabla^{n}}(x)}& = &
\displaystyle\frac{e^x_q}{[n]!}\, .
\end{array}
\end{equation}
Since $e^x_q$ is an increasing function, then $\displaystyle
\frac{f^{\nabla^{n}}(x)}{g^{\nabla^{n}}(x)}$ is increasing on
$[\sigma(a),b]$. From Theorem~\ref{hospital thm nabla} and the
relationships listed in (\ref{formulas}) we deduce that each of
the following functions
$$
\frac{f^{\nabla^{n}}(x)}{g^{\nabla^{n}}(x)}, \displaystyle
\frac{f^{\nabla^{n-1}}(x)}{g^{\nabla^{n-1}}(x)}, \ldots,
\frac{f^{\nabla}(x)}{g^{\nabla}(x)}, \frac{f(x)}{g(x)}
$$
is increasing on $[\sigma(a),b]$. Therefore,
for each $x \in [\sigma(a),b]$,
$$
\frac{f(\sigma(a))}{g(\sigma(a))} \leq \frac{f(x)}{g(x)} \leq
\frac{f(b)}{g(b)},
$$
proving the intended inequalities.
\end{proof}


\section*{Acknowledgment}

This work was partially supported by the R\&D unit
``Center for Research and Development in Mathematics and Applications''
(CIDMA) of the University of Aveiro,
cofinanced by the European Community Fund FEDER/POCI 2010.



\end{document}